\documentclass{article}%
\usepackage{amsmath}
\usepackage{amsfonts}
\usepackage{amssymb}
\usepackage{graphicx}%
\providecommand{\U}[1]{\protect\rule{.1in}{.1in}}

\begin{document}

\title{Relative belief inferences from decision theory}
\author{Michael\ Evans\\Department of Statistical Sciences, University of Toronto
\and Gun Ho Jang\\Ontario Institute for Cancer Research}
\date{}
\maketitle

\begin{abstract}
Relative belief inferences are shown to arise as Bayes rules or limiting Bayes
rules. These inferences are invariant under reparameterizations and possess a
number of optimal properties. In particular, relative belief inferences are
based on a direct measure of statistical evidence.

\end{abstract}


\section{Introduction}

Consider a sampling model for data $x$, given by a collection of densities
$\{f(\cdot\,|\,\theta):\theta\in\Theta\}$ with respect to a support measure
$\mu$ on sample space $\mathcal{X},$ and a proper prior, given by density
$\pi$ with respect to support measure $\nu$ on $\Theta.$ When the data
$x\in\mathcal{X}$ is observed these ingredients lead to the posterior
distribution on $\Theta$ with density given by $\pi(\theta\,|\,x)=\pi
(\theta)f(x\,|\,\theta)/m(x)$ with respect to support measure $\nu$ where
$m(x)=\int_{\Theta}\pi(\theta)f(x\,|\,\theta)\,\nu(d\theta)$ is the prior
predictive density of the data. In addition, there is a quantity of interest
$\psi=\Psi(\theta),$ where $\Psi:\Theta\rightarrow\Psi(\Theta)$ for which
inferences, such as an estimate $\psi(x)$ or a hypothesis assessment
$H_{0}:\Psi(\theta)=\psi_{0},$ are required. Let $\pi_{\Psi}$ denote the
marginal prior density of $\psi$ and $m(x\,|\,\psi)=\int_{\Theta
}f(x\,|\,\theta)\,\Pi(d\theta\,|\,\psi)$ be the conditional prior predictive
of the data after integrating out the nuisance parameters via the conditional
distribution of $\theta$ given $\Psi(\theta)=\psi.$ Bayesian inferences for
$\psi$ are then based on the ingredients $(\{m(\cdot\,|\,\psi):\psi\in
\Psi(\Theta)\},\pi_{\Psi},x)$\ alone or by adding a loss function $L.$

There are several different general approaches to deriving inferences based on
such ingredients. The two most commonly used are MAP-based inferences and
Bayesian decision theory. MAP-based inferences are based implicitly on
assuming that posterior probabilities can measure statistical evidence and do
not use a loss function explicitly. Bayesian decision theory seeks inferences
that are optimal with respect to risk which is defined as the expected loss
incurred by an inference under the joint distribution of $(\theta,x)\sim
\pi(\theta)f(x\,|\,\theta).$ Such optimal inferences are referred to as Bayes
rules and they generally exist. A concern with MAP-based inferences is that it
is not clear that posterior probabilities do measure evidence in addition to
measuring belief. A concern with decision-theory inferences is that, while the
model and prior are checkable against the data through model checking and
checking for prior-data conflict, it is not clear how to check the loss
function which can be viewed as being a somewhat arbitrary choice. In both
cases this renders such inferences of questionable validity for scientific applications.

Another approach to deriving Bayesian inferences is through relative belief.
Relative belief refers to how beliefs change from a priori to a posteriori.
This leads to a more natural approach to characterizing statistical evidence:
since it is the data that leads to change in beliefs from a priori to a
posteriori, it is this change that tells us whether evidence has been found in
favor of or against some specific value $\psi.$ In essence, in this approach
it is the a posteriori beliefs relative to the a priori beliefs that determine
inferences and not just a posteriori beliefs alone. Also, a loss function
plays no role in determining the inferences. Using relative belief as the
basis for deriving inferences produces statistical methodology with a number
of attractive features.

A historical theme in statistical research has been to seek an acceptable
definition of statistical evidence and, once found, use this to derive
inferences. For example, this is the focus of much of the work of Alan
Birnbaum, see Birnbaum (1962), who sought such a definition within the context
of frequentist inference. Frequentist theory, as opposed to Bayesian theory,
uses the ingredients $(\{f(\cdot\,|\,\theta):\theta\in\Theta\},x)$ together
with the idea that inferences should be graded according to their behavior in
hypothetical repeated sampling experiments and, hopefully this would lead to a
prescription of the inferences. Despite some impressive accomplishments, it is
fair to say that Birnbaum's program did not succeed as there is still no such
generally acceptable definition of statistical evidence within the frequentist
context. The pure likelihood theory of Royall (1997) is also concerned with
basing inference on a definition of statistical evidence, and also uses just
the ingredients $(\{f(\cdot\,|\,\theta):\theta\in\Theta\},x).$ Pure likelihood
theory invokes the likelihood principle to assert that the likelihood function
itself is the appropriate characterization of statistical evidence and bases
all inferences on the likelihood with no appeal to repeated sampling.
Frequency theory and likelihood theory have some appealing characteristics,
but both leave gaps in their approach to statistical reasoning. In particular,
inferences for marginal parameters $\psi=\Psi(\theta)$ can be problematical.
These issues are discussed in more detail in Evans (2024).

While not concerned directly with statistical evidence, Bayesian decision
theory has some obvious virtues. In particular, there is an axiomatization due
to Savage (1971). These axioms suggest that to not follow this path in
carrying out a statistical analysis is to commit an error. One may not find
the specific axioms of Savage acceptable, but it is difficult to argue that
the subject of Statistics does not need such an axiomatic formulation as
otherwise almost any statistical analysis seems justifiable.

The addition of a prior is what leads to a clear definition of statistical
evidence and so, provided one accepts the usage of priors, relative belief
essentially solves Birnbaum's problem. The purpose of this paper is to review
and extend results that show that relative belief inferences can be considered
as arising within the context of Bayesian decision theory. This, of course
requires the use of a loss function and it will be seen that the loss
functions that are used are checkable against the data and so appropriate for
scientific applications. So these inferences satisfy two of the great themes
of statistical research over the years, namely, they are evidence based and
yet justifiable within the context of decision theory.

It can be shown, for example, see Bernardo and Smith (2000), that MAP
inferences arise as the limits of Bayes rules via a sequence of loss functions%
\begin{equation}
L_{\lambda}(\theta,\psi)=I_{B_{_{\lambda}}^{c}(\psi)}(\Psi(\theta))\label{eq4}%
\end{equation}
where $\lambda>0$ and $B_{\lambda}(\Psi(\theta))$ is the ball of radius
$\lambda$ centered at $\psi.$ These inferences, however, are not invariant
under reparameterizations. It is shown here that relative belief inferences
also arise via a sequence of loss functions similar to (\ref{eq4}) but based
on the prior and these inferences are invariant. In general, Bayes rules will
also not be invariant under reparameterizations. Robert (1996) proposed using
the intrinsic loss function based on a measure of distance between sampling
distributions as Bayes rules with respect to such losses are invariant.
Bernardo (2005) proposed using the intrinsic loss function based on the
Kullback-Leibler divergence $KL(f_{\theta},f_{\theta^{\prime}})$ between
$f(\cdot\,|\,\theta)$ and $f(\cdot\,|\,\theta^{\prime}).$ When $\psi=\theta$
the intrinsic loss function is given by $L(\theta,\theta^{\prime}%
)=\min(KL(f(\cdot\,|\,\theta),f(\cdot\,|\,\theta^{\prime})),KL(f(\cdot
\,|\,\theta^{\prime}),f(\cdot\,|\,\theta))).$ For a general marginal parameter
$\psi$ the intrinsic loss function is $L(\theta,\psi)=\inf_{\theta^{\prime}%
\in\Psi^{-1}\{\psi\}}L(\theta,\theta^{\prime}).$ These loss functions are
intrinsic because they are based on the sampling model alone and they are
checkable via model checking.

Another possibility for an intrinsic loss function is to base the loss
function on the prior and this is in essence how the loss function arises in
the relative belief context. It is to be stressed, however, that the essential
ingredient of this approach is the clear characterization of what is meant by
statistical evidence and the loss function approach is not essential for its
justification. It is, however, a satisfying result that relative belief can be
placed into the decision-theoretic context with the loss being checkable via
checking for prior-data conflict. Furthermore, the loss function used has some
direct appeal.

In some contexts relative belief inferences are Bayes rules, but in a general
context they are seen to arise as the limits of Bayes rules. This approach has
some historical antecedents. For example, in Le Cam (1953) it is shown that
the MLE is asymptotically Bayes but this is for a fixed loss function, with
increasing amounts of data and a sequence of priors. In the context discussed
here there is a fixed amount of data, a fixed model and prior but there is a
sequence of loss functions all based on the single fixed prior.

While it is preferable in many applications to state the inferences solely
based on the evidence in the data, one can still consider inferences that
possess some kind of optimality with respect to loss. Any discrepancy can then
be justified based on particular characteristics of the application, e.g.,
evidence is obtained that a drug generally prevents the progression of a
disease but the expense and side effects are too great to warrant its usage.
So it is not being suggested here that decision-theoretic inferences are not
relevant as indeed the concept of utility or loss is a significant component
of many applications.

Section 2 is concerned with describing the general characteristics of the
three approaches to deriving Bayesian inferences. Section 3 shows how relative
belief estimation and prediction inferences can be seen to arise from decision
theory and Section 4 does this for credible regions and hypothesis assessment.
Throughout the paper the probability measures associated with a density are
denoted by the same letter but capitalized. All proofs of theorems and
corollaries are in the Appendix excepting the case where $\Psi(\Theta)$ is
finite as these are quite straightforward and supply motivation for the more
complicated contexts. The overall goal of the paper is to show that relative
belief inferences can arise through decision-theoretic considerations even
though their primary motivation is through characterizing statistical
evidence. In particular, it is shown here that relative belief estimators, as
used in practice, are admissible. Some of the discussion here has appeared in
the book Evans (2015) and is included to provide a complete exposition of this relationship.

\section{Bayesian inference}

The various approaches to deriving Bayesian inferences are now described in
some detail.

\subsection{MAP inferences}

The highest posterior density (hpd), or MAP-based, approach to determining
inferences constructs credible regions of the form
\begin{equation}
H_{\gamma}(x)=\{\psi:\pi_{\Psi}(\psi\,|\,x)\geq h_{\gamma}(x)\} \label{eq1}%
\end{equation}
where $\pi_{\Psi}(\cdot\,|\,x)$ is the marginal posterior density with respect
to a support measure $\nu_{\Psi}$ on $\Psi(\Theta),$ and $h_{\gamma}(x)$ is
chosen so that $h_{\gamma}(x)=\sup\{k:\Pi_{\Psi}(\{\psi:\pi_{\Psi}%
(\psi\,|\,x)\geq k\}\,|\,x)\geq\gamma\}.$ It follows from (\ref{eq1}) that, to
assess the hypothesis $H_{0}:\Psi(\theta)=\psi_{0},$ then we can use the tail
probability given by $1-\inf\{\gamma:\psi_{0}\in H_{\gamma}(x)\}.$
Furthermore, the class of sets $H_{\gamma}(x)$ is naturally "centered" at the
posterior mode (when it exists uniquely) as $H_{\gamma}(x)$ converges to this
point as $\gamma\rightarrow0.$ The use of the posterior mode as an estimator
is commonly referred to as MAP (maximum \textit{a posteriori}) estimation. We
can then think of the size of the set $H_{\gamma}(x),$ say for $\gamma=0.95,$
as a measure of how accurate the MAP estimator is in a given context.
Furthermore, when $\Psi(\Theta)$ is an open subset of a Euclidean space, then
$H_{\gamma}(x)$ minimizes volume among all $\gamma$-credible regions.

It is well-known, however, that hpd inferences suffer from a defect. In
particular, in the continuous case MAP inferences are not invariant under
reparameterizations. For example, this means that if $\psi_{MAP}(x)$ is the
MAP estimate of $\psi$, then it is not necessarily true that $\Upsilon
(\psi_{MAP}(x))$ is the MAP estimate of $\tau=\Upsilon(\psi)$ when $\Upsilon$
is a 1-1, smooth transformation. The noninvariance of a statistical procedure
seems very unnatural as it implies that the statistical analysis depends on
the parameterization and typically there does not seem to be a good reason for
this. Note too that estimates based upon taking posterior expectations will
also suffer from this lack of invariance. It is also the case that
MAP\ inferences are not based on a direct characterization of statistical
evidence. Both of these issues motivate the development of relative belief inferences.

\subsection{Bayesian decision theory}

An ingredient that is commonly added to $(\{f(\cdot\,|\,\theta):\theta
\in\Theta\},\pi,x)$ is a loss function, namely, $L:\Theta\times\Psi
(\Theta)\rightarrow\lbrack0,\infty)$ satisfying $L(\theta,\psi)=L(\theta
^{\prime},\psi)$ whenever $\Psi(\theta)=\Psi(\theta^{\prime})$ and
$L(\theta,\psi)=0$ only when $\psi=\Psi(\theta).$ The goal is to find a
procedure, say $\delta(x)\in\Psi(\Theta),$ that in some sense minimizes the
loss $L(\theta,\delta(x))$ based on the joint distribution of $(\theta,x).$
Given the assumptions on the loss function the loss function can instead be
thought of as a map $L:\Psi(\Theta)\times\Psi(\Theta)\rightarrow
\lbrack0,\infty)$ with $L(\psi,\psi^{\prime})=0\,$iff $\psi=\psi^{\prime}$ and
the ingredients can be represented as $(\{m(\cdot\,|\,\psi):\psi\in\Psi
(\Theta)\},\pi_{\Psi},L,x).$

The goal of a decision analysis is then to find a decision function
$\delta:\mathcal{X}\rightarrow\Psi(\Theta)$ that minimizes the \textit{prior
risk}
\[
r(\delta)=\int_{\Psi(\Theta)}\int_{\mathcal{X}}L(\psi,\delta
(x))\,M(dx\,|\,\psi)\,\Pi_{\Psi}(d\psi)=\int_{\mathcal{X}}r(\delta
\,|\,x)\,M(dx)\,
\]
where $r(\delta\,|\,x)=\int_{\Psi(\Theta)}L(\psi,\delta(x))\,\Pi_{\Psi}%
(d\psi\,|\,x)$ is the \textit{posterior risk}. Such a $\delta$ is called a
Bayes rule and clearly a $\delta$ that minimizes $r(\delta\,|\,x)$ for each
$x$ is a Bayes rule. Further discussion of Bayesian decision theory can be
found in Berger (1985).

As noted in Bernardo (2005) a decision formulation also leads to credible
regions for $\psi,$ namely, a $\gamma$\textit{-lowest posterior loss credible
region} is defined by
\begin{equation}
D_{\gamma}(x)=\{\psi:r(\psi\,|\,x)\leq d_{\gamma}(x)\} \label{eq3}%
\end{equation}
where $d_{\gamma}(x)=\inf\{k:\Pi_{\Psi}(\{\psi:r(\psi\,|\,x)\leq
k\}\,|\,x)\geq\gamma.$ Note that $\psi$ in (\ref{eq3}) is interpreted as the
decision function that takes the value $\psi$ constantly in $x.$ Clearly as
$\gamma\rightarrow0$ the set $D_{\gamma}(x)$ converges to the value of a Bayes
rule at $x.$ For example, with quadratic loss the Bayes rule is given by the
posterior mean and a $\gamma$-lowest posterior loss region is the smallest
sphere centered at the mean containing at least $\gamma$ of the posterior probability.

\subsection{Relative belief inferences}

Relative belief inferences, like MAP inferences, are based on the ingredients
$(\{m(\cdot\,|\,\psi):\psi\in\Psi(\Theta)\},\pi_{\Psi},x).$ Note that
underlying both approaches is the principle (axiom) of conditional probability
that says that initial beliefs about $\psi,$ as expressed by the prior
$\pi_{\Psi},$ must be replaced by conditional beliefs as expressed by he
posterior $\pi_{\Psi}(\cdot\,|\,x).$ In this approach, however, a measure of
statistical evidence is used given by the relative belief ratio%
\begin{equation}
RB_{\Psi}(\psi\,|\,x)=\frac{\pi_{\Psi}(\psi\,|\,x)}{\pi_{\Psi}(\psi)}%
=\frac{m(x\,|\,\psi)}{m(x)}. \label{rb1}%
\end{equation}
The relative belief ratio produces the following conclusions: if $RB_{\Psi
}(\psi\,|\,x)>1,$ then there is evidence in favor of $\psi$ being the true
value, if $RB_{\Psi}(\psi\,|\,x)<1,$ there is evidence against $\psi$ being
the true value and if $RB_{\Psi}(\psi\,|\,x)=1,$ then there is no evidence
either way. These implications follow from a very simple principle of inference.

\begin{quote}
Principle of Evidence: for probability model $(\Omega,\mathcal{F},P),$ if
$C\in\mathcal{F}$ is observed to be true where $P(C)>0,$ then there is
evidence in favor of $A\in\mathcal{F}$ being true if $P(A\,|\,C)>P(A),$
evidence against $A\in\mathcal{F}$ being true if $P(A\,|\,C)<P(A)$ and no
evidence either way if $P(A\,|\,C)=P(A).$
\end{quote}

\noindent This principle seems obvious when $\Pi_{\Psi}$ is a discrete
probability measure. For the continuous case, where $\Pi_{\Psi}(\{\psi\})=0,$
let $N_{\epsilon}(\psi)$ be a sequence of neighborhoods of $\psi$\ converging
nicely to $\psi$ as $\epsilon\rightarrow0$ (see Rudin (1974)), then under weak
conditions, e.g., $\pi_{\Psi}$ is continuous and positive at $\psi,$
\[
\lim_{\epsilon\rightarrow0}RB_{\Psi}(N_{\epsilon}(\psi)\,|\,x)=\lim
_{\epsilon\rightarrow0}\frac{\Pi_{\Psi}(N_{\epsilon}(\psi)\,|\,x)}{\Pi_{\Psi
}(N_{\epsilon}(\psi))}=\frac{\pi_{\Psi}(\psi\,|\,x)}{\pi_{\Psi}(\psi
)}=RB_{\Psi}(\psi\,|\,x)
\]
and this justifies the general interpretation of $RB_{\Psi}(\psi\,|\,x)$ as a
measure of evidence. The relative belief ratio determines the inferences.

A natural estimate of $\psi$ is the \textit{relative belief estimate}
\[
\psi_{RB}(x)=\arg\sup_{\psi}RB_{\Psi}(\psi\,|\,x)
\]
as it has the maximum evidence in favor. To assess the accuracy of $\psi
_{RB}(x)$ there is the \textit{plausible region} $Pl_{\Psi}(x)=\{\psi
:RB_{\Psi}(\psi\,|\,x)>1\},$ the set of $\psi$ values having evidence in favor
of being the true value. The size of $Pl_{\Psi}(x)$ together with its
posterior content $\Pi_{\Psi}(Pl_{\Psi}(x)\,|\,x),$ which measures the belief
that the true value is in $Pl_{\Psi}(x),$ provide the assessment of the
accuracy$.$ So, if $Pl_{\Psi}(x)$ is "small" and $\Pi_{\Psi}(Pl_{\Psi
}(x)\,|\,x)\approx1,$ then $\psi_{RB}(x)$ is to be considered as an accurate
estimate of $\psi$ but not otherwise. A relative belief $\gamma$-credible
region
\[
C_{\Psi,\gamma}(x)=\{\psi:RB_{\Psi}(\psi\,|\,x)\geq c_{\gamma}(x)\},
\]
where $c_{\gamma}(x)=\sup\{c:\Pi_{\Psi}(RB_{\Psi}(\psi\,|\,x)\geq
c\,|\,x)\geq\gamma\},$ for $\psi$ can also be quoted provided $C_{\Psi,\gamma
}(x)\subset Pl_{\Psi}(x).$ The containment is necessary as otherwise
$C_{\Psi,\gamma}(x)$ would contain a value $\psi$ for which there is evidence
against $\psi$ being the true value.

For assessing the hypothesis $H_{0}:\Psi(\theta)=\psi_{0},$ the value
$RB_{\Psi}(\psi_{0}\,|\,x)$ indicates whether there is evidence in favor of or
against $H_{0}.$ The strength of this evidence can be measured by the
posterior probability $\Pi_{\Psi}(\{\psi_{0}\}\,|\,x),$ as this measures the
belief in what the evidence says. So, if $RB_{\Psi}(\psi_{0}\,|\,x)>1$ and
$\Pi_{\Psi}(\{\psi_{0}\}\,|\,x)\approx1,$ then there is strong evidence that
$H_{0}$ is true while, when $RB_{\Psi}(\psi_{0}\,|\,x)<1$ and $\Pi_{\Psi
}(\{\psi_{0}\}\,|\,x)\approx0,$ there is strong evidence that $H_{0}$ is
false. Since $\Pi_{\Psi}(\{\psi_{0}\})$ can be small, even 0 in the continuous
case, it makes more sense to measure the strength of the evidence in such a
case by
\[
Str_{\Psi}(\psi_{0}\,|\,x)=\Pi_{\Psi}(RB_{\Psi}(\psi\,|\,x)\leq RB_{\Psi}%
(\psi_{0}\,|\,x)\,|\,x).
\]
If $RB_{\Psi}(\psi_{0}\,|\,x)>1$ and $Str_{\Psi}(\psi_{0}\,|\,x)\approx1,$
then the evidence is strong that $\psi_{0}$ is the true values as there is
small belief that the true value of $\psi$ has more evidence in its favor than
$\psi_{0}.$ If $RB_{\Psi}(\psi_{0}\,|\,x)<1$ and $Str_{\Psi}(\psi
_{0}\,|\,x)\approx0,$ then the evidence is strong that $\psi_{0}$ is not the
true values as there is large belief that the true value of $\psi$ has more
evidence in its favor than $\psi_{0}.$ Actually, there is no reason to quote a
single number to measure the strength and both $\Pi_{\Psi}(\{\psi
_{0}\}\,|\,x)$ and $Str_{\Psi}(\psi_{0}\,|\,x)$ can be quoted when relevant.

An important aspect of both $Str_{\Psi}(\psi_{0}\,|\,x)$ and $\Pi_{\Psi
}(Pl_{\Psi}(x)\,|\,x)$ is what happens as the amount of data increases. To
ensure that these behave appropriately, namely, $Str_{\Psi}(\psi
_{0}\,|\,x)\rightarrow0(1)$ when $H_{0}$ is false(true) and $\Pi_{\Psi
}(Pl_{\Psi}(x)\,|\,x)\rightarrow1,$ it is necessary to take into account the
difference that matters $\delta.$ By this we mean that there is a distance
measure $d_{\Psi}$ on $\Psi(\Theta)\times\Psi(\Theta)$ such that if $d_{\Psi
}(\psi,\psi^{\prime})\leq\delta,$ then in terms of the application, these
values are considered equivalent. Such a $\delta$ always exists because
measurements are always taken to finite accuracy. For example, if $\psi$ is
real-valued, then there is a grid of values $\ldots\psi_{-2},\psi_{-1}%
,\psi_{0},\psi_{1},\psi_{2},\ldots$ separated by $\delta$ and inferences are
determined using the relative belief ratios of the intervals $[\psi_{i}%
-\delta/2,\psi_{i}+\delta/2).$ In effect, $H_{0}$ is now $H_{0}:\Psi
(\theta)\in\lbrack\psi_{0}-\delta/2,\psi_{0}+\delta/2).$ When the computations
are carried out in this way then $Str_{\Psi}(\psi_{0}\,|\,x)$ and $\Pi_{\Psi
}(Pl_{\Psi}(x)\,|\,x)$ do what is required. As a particular instance of this
see the results in Section 4 where such a discretization plays a key role.

It is easy to see that the class of relative belief credible regions
$\{C_{\Psi,\gamma}(x):\gamma\in\lbrack0,1]\}$ for $\psi$ is independent of the
marginal prior $\pi_{\Psi}.$ When a value $\gamma\in\lbrack0,1]$ is specified,
however, the set $C_{\Psi,\gamma}(x)$ does depend on $\pi_{\Psi}$ through
$c_{\gamma}(x).$ So the form of relative belief inferences about $\psi$ is
completely robust to the choice of $\pi_{\Psi}$ but the quantification of the
uncertainty in the inferences is not. For example, when $\psi=\Psi
(\theta)=\theta,$ then $\theta_{RB}(x)$ is the MLE while, in general,
$\psi_{RB}(x)$ is the maximizer of the integrated likelihood $m(x\,|\,\psi).$
Similarly, relative belief regions are likelihood regions in the case of the
full parameter, and integrated likelihood regions generally. As such,
likelihood regions can be seen as essentially Bayesian in character with a
clear and precise characterization of evidence through the relative belief
ratio and now have probability assignments through the posterior. It is the
case, however, that a relative belief ratio $RB_{\Psi}(\psi\,|\,x),$ while
proportional to an integrated likelihood, cannot be multiplied by an arbitrary
positive constant, as with a likelihood, without losing its interpretation in
measuring statistical evidence. It has been established in Al Labadi and Evans
(2017) that relative belief inferences for $\psi$ are optimally robust to the
prior $\pi_{\Psi}.$

As can be seen from (\ref{rb1}), relative belief inferences are always
invariant under smooth reparameterizations and this is at least one reason why
they are preferable to MAP inferences. It is the case, however, that any rule
for measuring evidence which satisfies the principle of evidence also produces
valid estimates as these lie in $Pl_{\Psi}(x)$ and so will have the same
"accuracy" as $\psi_{RB}(x).$ For example, if instead of the relative belief
ratio the difference $\pi_{\Psi}(\psi\,|\,x)-\pi_{\Psi}(\psi)$ was used as the
measure of evidence with cut-off 0, then this satisfies the principle of
evidence but the estimate is no longer necessarily invariant under
reparameterizations. The Bayes factor with cut-off 1 is also a valid measure
of evidence but there are a number of reasons why the relative belief ratio is
to be preferred to the Bayes factor for general inferences, see Al-Labadi,
Alzaatreh and Evans (2024).

\section{Estimation: discrete parameter space}

\label{sec:prior-loss}

The following theorem presents the basic definition of the loss function when
$\Psi$ is finite and establishes an important optimality result. The indicator
function for the set $A$ is denoted $I_{A}.$\smallskip

\noindent\textbf{Theorem 1}. Suppose that $\pi_{\Psi}(\psi)>0$ for every
$\psi\in\Psi(\Theta)$ where $\Psi(\Theta)$ is finite with $\nu_{\Psi}$ equal
to counting measure on $\Psi(\Theta).$ Then for the loss function%
\begin{equation}
L_{RB}(\theta,\psi)=\frac{I_{\{\psi\}^{c}}(\Psi(\theta))}{\pi_{\Psi}%
(\Psi(\theta))}, \label{eq5}%
\end{equation}
the relative belief estimator $\psi_{RB}$\ is a Bayes rule.

\noindent Proof: We have that\
\begin{align}
r(\delta\,|\,x)  &  =\int_{\Psi(\Theta)}\frac{I_{\{\delta(x)\}^{c}}(\psi)}%
{\pi_{\Psi}(\psi)}\,\Pi_{\Psi}(d\psi\,|\,x)\nonumber\\
&  =\int_{\Psi(\Theta)}RB_{\Psi}(\psi\,|\,x)\,\nu_{\Psi}(d\psi)-RB_{\Psi
}(\delta(x)\,|\,x). \label{eq5a}%
\end{align}
Since $\Psi(\Theta)$ is finite, the first term in (\ref{eq5a}) is finite and a
Bayes rule at $x$ is given by the value $\delta(x)$ that maximizes the second
term.\ Therefore, $\psi_{RB}$ is a Bayes rule. $\blacksquare$\smallskip

The loss function $L_{RB}$ seems very natural. For beliefs about the true
value of $\psi$ are expressed by the prior $\pi_{\Psi}$ and so values where
$\pi_{\Psi}(\psi)$ is very low and $\psi$ is indeed a false value, would be
quite misleading if the inferences pointed to such a value. So it is
appropriate for such values to bear large losses. In a sense the statistician
is acknowledging what such values are by the choice of prior. Of course, the
prior may be wrong in the sense that the bulk of its mass is placed in a
region where the true value of $\psi$ does not lie. This is why checking for
prior-data conflict, before inference is carried out, is always recommended.
Procedures for checking a prior are discussed in Evans and Moshonov (2006) and
Nott et al. (2020) and an approach to replacing a prior found to be at fault
is developed in Evans and Jang (2011). The loss $L_{RB}$ motivates the other
losses for relative belief discussed here so this comment applies to those
losses as well.

The prior risk of $\delta$ satisfies
\begin{equation}
r(\delta)=\int_{\Psi(\Theta)}\int_{\mathcal{X}}\frac{I_{\{\delta(x)\}^{c}%
}(\psi)}{\pi_{\Psi}(\psi)}\,M(dx\,|\,\psi)\,\Pi_{\Psi}(d\psi)=\sum_{\psi
\in\Psi(\Theta)}M(\delta(x)\neq\psi\,|\,\psi), \label{eq5b}%
\end{equation}
the sum of the conditional prior error probabilities over all $\psi$ values.
If instead the loss function is taken to be $L_{MAP}(\theta,\psi
)=I_{\{\psi\}^{c}}(\Psi(\theta)),$ as in (\ref{eq4}), then virtually the same
proof as Theorem 1 establishes that $\psi_{MAP}$ is a Bayes rule with respect
to this loss and the prior risk equals%
\begin{equation}
\sum_{\psi}M(\delta(x)\neq\psi\,|\,\psi)\pi_{\Psi}(\psi), \label{eq5c}%
\end{equation}
the prior probability of making an error. Both $L_{MAP}$ and $L_{RB}$ are
two-valued loss functions but, when an incorrect decision is made, the loss is
constant in $\Psi(\theta)$ for $L_{MAP}$ while it equals the reciprocal of the
prior probability of $\Psi(\theta)$ for $L_{RB}$. So $L_{RB}$ penalizes an
incorrect decision much more severely when the true value of $\Psi(\theta)$ is
in the tails of the prior. Note that $\psi_{MAP}=\psi_{RB}$ when $\Pi_{\Psi}$
is uniform. It is seen too that (\ref{eq5b}) is an upper bound on (\ref{eq5c})
so controlling losses based on $L_{RB}$ automatically controls the losses
based on $L_{MAP}.$

As already noted, $RB_{\Psi}(\psi\,|\,x)$ is proportional to the integrated
likelihood of $\psi.$ So, under the conditions of Theorem 1, the maximum
integrated likelihood estimator is a Bayes rule. Furthermore, the Bayes rule
is the same for every choice of $\pi_{\Psi}$ and only depends on the full
prior through the conditional prior $\Pi(\cdot\,|\,\psi)$ placed on the
nuisance parameters. When $\psi=\theta$ then $\theta_{RB}(x)$ is the MLE of
$\theta$ and so the MLE of $\theta$ is a Bayes rule for every prior $\pi.$

Note that when $\Psi(\Theta)=\{\psi_{0},\psi_{1}\}$ then $RB_{\Psi}(\psi
_{0}\,|\,x)>(<)1$ iff $RB_{\Psi}(\psi_{1}\,|\,x)<(>)1$ so $\psi_{RB}%
(x)=\psi_{0}$ when $RB_{\Psi}(\psi_{0}\,|\,x)>1$ and $\psi_{RB}(x)=\psi_{1}$
otherwise. This is the classical context for hypothesis testing and $\psi
_{RB}(x)=\psi_{0}$ can be viewed as acceptance of the hypothesis $H_{0}%
:\theta\in\Psi^{-1}\{\psi_{0}\}$ and $\psi_{RB}(x)=\psi_{1}$ as rejection of
$H_{0}.$ Theorem 1 establishes that relative belief provides a Bayes rule for
the hypothesis testing problem.

The loss function (\ref{eq5}) does not provide meaningful results when
$\Psi(\Theta)$ is infinite as (\ref{eq5b}) shows that $r(\delta)$ will be
infinite. So we modify (\ref{eq5}) via a parameter $\eta>0$ and define the
loss function
\begin{equation}
L_{RB,\eta}(\theta,\psi)=\frac{I_{\{\psi\}^{c}}(\Psi(\theta))}{\max(\eta
,\pi_{\Psi}(\Psi(\theta)))}. \label{eq6}%
\end{equation}
Note that $L_{\eta,RB}$ is a bounded by $1/\eta.$ This loss function is like
(\ref{eq5}) but does not allow for arbitrarily large losses. Without loss of
generality we can restrict $\eta$ to a sequence of values converging to
0.\smallskip

\noindent\textbf{Theorem 2}. Suppose that $\pi_{\Psi}(\psi)>0$ for every
$\psi\in\Psi(\Theta),$ that $\Psi(\Theta)$ is countable with $\nu_{\Psi}$
equal to counting measure and that $\psi_{RB}(x)$ is the unique maximizer of
$RB_{\Psi}(\psi\,|\,x)$ for all $x.$ For the loss function (\ref{eq6}) Bayes
rule $\delta_{\eta},$ then $\delta_{\eta}(x)\rightarrow\psi_{RB}(x)$ as
$\eta\rightarrow0,$ for every $x\in\mathcal{X}.$\smallskip

\noindent The proof of Theorem 2 also establishes the following
results.\smallskip

\noindent\textbf{Corollary 3}. For all sufficiently small $\eta$ the value of
a Bayes rule at $x$ is given by $\psi_{RB}(x).$\smallskip

The following is an immediate consequence of Theorem 1 and Corollary 3 as
$\psi_{RB}$ is a Bayes rule.\smallskip

\noindent\textbf{Corollary 4}. $\psi_{RB}$ is an admissible estimator with
respect to the loss $L_{RB}$ when $\Psi(\Theta)$ is finite and the loss
$L=L_{RB,\eta},$ with $\eta$ sufficiently small, when $\Psi(\Theta)$ is
countable. \smallskip

In a general estimation problem $\delta$ is risk unbiased with respect to a
loss function $L$ if $E_{\theta}(L(\theta^{\prime},\delta(x)))\geq E_{\theta
}(L(\theta,\delta(x)))$ for all $\theta^{\prime},\theta\in\Theta.$ This says
that on average $\delta(x)$ is closer to the true value than any other value
when we interpret $L(\theta,\delta(x))$ as a measure of distance between
$\delta(x)$ and $\Psi(\theta).$ A definition of \textit{Bayesian unbiasedness}
for $\delta$ with respect to \thinspace$L$ is that
\[
\int_{\Theta}\int_{\Theta}E_{\theta}(L(\theta^{\prime},\delta(x)))\,\Pi
(d\theta)\,\Pi(d\theta^{\prime})\geq\int_{\Theta}E_{\theta}(L(\theta
,\delta(x)))\,\Pi(d\theta)=r(\delta)
\]
as this retains the idea of being closer on average to the true value than a
false value. Consider now a family of loss functions of the form
\begin{equation}
L(\theta,\psi)=I_{\{\psi\}^{c}}(\Psi(\theta))h(\Psi(\theta)) \label{eq6c}%
\end{equation}
where $h$ is a nonnegative function satisfying $\int_{\Theta}h(\Psi
(\theta))\,\Pi(d\theta)<\infty\ $and\thinspace\ note that this includes
$L_{RB}$ and $L_{MAP}$ when $\Psi(\Theta)$ is finite and $L_{RB,\eta}%
.$\smallskip

\noindent\textbf{Theorem 5}. If $\Psi(\Theta)$ is finite or countable, then
$\psi_{RB}(x)$ is Bayesian unbiased under the loss function (\ref{eq6c}%
).\smallskip

Suppose after observing $x$ it is desired to predict a future (or concealed)
value $y\in\mathcal{Y}$ where $y\sim g_{\delta(\theta)}(y\,|\,x),$ a density
with respect to support measure $\mu_{\mathcal{Y}}$ on $\mathcal{Y},$ and it
is assumed that the true value of $\theta$ in the model for $x$ gives the true
value of $\delta(\theta).$ The prior predictive density of $y$ is given by
$q(y)=\int_{\Theta}\int_{\mathcal{X}}\pi(\theta)f_{\theta}(x)g_{\delta
(\theta)}(y\,|\,x)\,\mu(dx)\,\nu(d\theta)$ while the posterior predictive
density is $q(y\,|\,x)=\int_{\Theta}\pi(\theta\,|\,x)g_{\delta(\theta
)}(y\,|\,x)\,\nu(d\theta).$ The relative belief ratio for a future value $y$
is thus $RB_{Y}(y\,|\,x)=q(y\,|\,x)/q(y)$ and the relative belief prediction,
namely, the value maximizing $RB_{Y}(\cdot\,|\,x),$ is denoted $y_{RB}(x).$
When $\mathcal{Y}$ is finite then, with basically the same argument as in
Theorem 1, $y_{RB}$ is a Bayes rule under the loss function $L_{RB}%
(y,y^{\prime})=I_{\{y\}^{c}}(y^{\prime})/q(y).$ Also, it can be proved that
$y_{RB}$ is a limit of Bayes rules when $\mathcal{Y}$ is countable.\smallskip

Consider now a common application where $\Psi(\Theta)$ is finite.\smallskip

\noindent\textbf{Example 1.} \textit{Classification}

For a classification problem there are $k$ categories $\{\psi_{1},\ldots
,\psi_{k}\},$ prescribed by some function $\Psi,$ where $\pi_{\Psi}(\psi
_{i})>0$ for each $i.$ Estimating $\psi$ is then equivalent to classifying the
data as having come from one of the distributions in the classes specified by
$\Psi^{-1}\{\psi_{i}\}.$ The standard Bayesian solution to this problem is to
use $\psi_{MAP}(x)$ as the classifier. From (\ref{eq5c}) we have that
$\psi_{MAP}(x)$ minimizes the prior probability of misclassification while
from (\ref{eq5b}) $\psi_{RB}(x)$ minimizes the sum of the probabilities of
misclassification. The essence of the difference is that $\psi_{RB}(x)$ treats
the errors of misclassification equally while $\psi_{MAP}(x)$ weights the
errors by their prior probabilities.

The following shows that minimizing the sum of the error probabilities is
often more appropriate than minimizing the weighted sum. Suppose $k=2$ and
$x\sim$ Bernoulli$(\psi_{0})$ or $x\sim$ Bernoulli$(\psi_{1})$ with $\pi
(\psi_{0})=1-\epsilon$ and $\pi(\psi_{1})=\epsilon$ representing the known
proportions of individuals either labelled coming from population 0 or 1. For
example, consider $\psi_{0}$ as the probability of a positive diagnostic test
for a disease in the nondiseased population while $\psi_{1}$ is this
probability for the diseased population. Suppose that $\psi_{0}/\psi_{1}$ is
very small, indicating that the test is successful in identifying the disease
while not yielding many false positives, and that $\epsilon$ is very small, so
the disease is rare. The question then is to assign a randomly chosen
individual to a population based on the results of their test.

The posterior is given by $\pi(\psi_{0}\,|\,1)=\psi_{0}(1-\epsilon)/(\psi
_{0}(1-\epsilon)+\psi_{1}\epsilon)$ and $\pi(\psi_{0}\,|\,0)=(1-\psi
_{0})(1-\epsilon)/((1-\psi_{0})(1-\epsilon)+(1-\psi_{1})\epsilon).$
Therefore,
\begin{align*}
\psi_{MAP}(1)  &  =\left\{
\begin{tabular}
[c]{ll}%
$\psi_{0}$ & if $\psi_{0}/\psi_{1}>\epsilon/(1-\epsilon)$\\
$\psi_{1}$ & otherwise
\end{tabular}
\ \right. \\
\psi_{MAP}(0)  &  =\left\{
\begin{tabular}
[c]{ll}%
$\psi_{0}$ & if $(1-\psi_{0})/(1-\psi_{1})>\epsilon/(1-\epsilon)$\\
$\psi_{1}$ & otherwise
\end{tabular}
\ \right.
\end{align*}
This implies that $\psi_{MAP}$ will always classify a person to the
nondiseased population when $\epsilon$ is small enough, e.g., when $\psi
_{0}=0.05,\psi_{1}=0.80,$ and $\epsilon<0.0625.$ By contrast, in this
situation, $\psi_{RB}$ always classifies an individual with a positive test to
the diseased population and to the nondiseased population for a negative test.
Since $M(\cdot\,|\,\psi_{i})$ is the Bernoulli$(\psi_{i})$ distribution, when
$\psi_{0}<\psi_{1}$ and $\epsilon$ is small enough,
\begin{align*}
M(\psi_{MAP}  &  \neq\psi_{0}\,|\,\psi_{0})+M(\psi_{MAP}\neq\psi_{1}%
\,|\,\psi_{1})=0+1=1,\\
M(\psi_{RB}  &  \neq\psi_{0}\,|\,\psi_{0})+M(\psi_{RB}\neq\psi_{1}%
\,|\,\psi_{1})=\psi_{0}+(1-\psi_{1})=0.25.
\end{align*}
This illustrates clearly the difference between these two procedures as
$\psi_{RB}$ does better than $\psi_{MAP}$ on the diseased population when
$\psi_{0}$ is small and $\psi_{1}$ is large as would be the case for a good
diagnostic. Of course $\psi_{MAP}$ minimizes the overall error rate but at the
price of ignoring the most important class in this problem. Note that this
example can be extended to the situation where we need to estimate the
$\psi_{i}$ based on samples from the respective populations but this will not
materially affect the overall conclusions.

Consider now a situation where $(x,c)$ is such that $x\,|\,c\sim
f_{c},c\,|\,\epsilon\sim$ Bernoulli$(\epsilon)$ where $f_{0}$ and $f_{1}$ are
known but $\epsilon$ is unknown with prior $\pi.$ This is a generalization of
the previous discussion where $\epsilon$ was assumed to be known. Then based
on a sample $(x_{1},c_{1}),\ldots,(x_{n},c_{n})$ from the joint distribution
the goal is to predict the value $c_{n+1}$ for a newly observed $x_{n+1}.$

The prior of $c$ is $q(c)=\int_{0}^{1}(1-\epsilon)^{1-c}\epsilon^{c}%
\pi(\epsilon)\,d\epsilon\ $and, if $\epsilon\sim$ beta$(\alpha,\beta),$ so the
prior predictive of $c_{n+1}$ is Bernoulli$(\alpha/(\alpha+\beta)).$ The
posterior predictive density of $c_{n+1}$ equals, where $\bar{c}=n^{-1}%
\sum_{i=1}^{n}c_{i},$%
\begin{align*}
&  q(c\,|\,(x_{1},c_{1}),\ldots,(x_{n},c_{n}),x_{n+1})\\
&  \propto(f_{0}(x_{n+1}))^{1-c}(f_{1}(x_{n+1}))^{c}\int_{0}^{1}%
\epsilon^{n\bar{c}+c}(1-\epsilon)^{n(1-\bar{c})+(1-c)}\pi(\epsilon
)\,d\epsilon\\
&  =f_{c}(x_{n+1})\Gamma\left(  \alpha+n\bar{c}+c\right)  \Gamma
(\beta+n(1-\bar{c})+1-c).
\end{align*}
It follows that, suppressing the dependence on the data,
\begin{align}
c_{MAP}  &  =\left\{
\begin{array}
[c]{cl}%
1 & \text{if }\frac{f_{1}(x_{n+1})}{f_{0}(x_{n+1})}\frac{\left(  \alpha
+n\bar{c}\right)  }{\left(  \beta+n(1-\bar{c})\right)  }>1\\
0 & \text{otherwise,}%
\end{array}
\right. \nonumber\\
c_{RB}  &  =\left\{
\begin{array}
[c]{cl}%
1 & \text{if }\frac{f_{1}(x_{n+1})}{f_{0}(x_{n+1})}\frac{\beta\left(
\alpha+n\bar{c}\right)  }{\alpha\left(  \beta+n(1-\bar{c})\right)  }<1\\
0 & \text{otherwise}%
\end{array}
\right.  \label{predclass}%
\end{align}
Note that $c_{MAP}$ and $c_{RB}$ are identical whenever $\alpha=\beta.$

From these formulas it is apparent that a substantial difference will arise
between $c_{MAP}$ and $c_{RB}$ when one of $\alpha$ or $\beta$ is much bigger
than the other. As in Example 1 these correspond to situations where we
believe that $\epsilon$ or $1-\epsilon$ is very small. Suppose we take
$\alpha=1$ and let $\beta$ be relatively large, as this corresponds to knowing
\textit{a priori} that $\epsilon$ is very small. Then (\ref{predclass})
implies that $c_{MAP}\leq c_{RB}$ and so $c_{RB}=1$ whenever $c_{MAP}=1.$ A
similar conclusion arises when we take $\beta=1$ and $\alpha<1.$

To see what kind of improvement is possible consider a simulation study. Let
$f_{0}$ be a $N(0,1)$ density, $f_{1}$ be a $N(\mu,1)$ density, $n=10$ and the
prior on $\epsilon$ be beta$(1,\beta).$ Table 1 presents the Bayes risks for
$c_{MAP}\ $and $c_{RB}$ for various choices of $\beta$ when $\mu=1.$ When
$\beta=1$ they are equivalent but we see that as $\beta$ rises the performance
of $c_{MAP}$ deteriorates while $c_{RB}$ improves. Large values of $\beta$
correspond to having information that $\epsilon$ is small. When $\beta=14$
about $0.50$ of the prior probability is to the left of $0.05,$ with
$\beta=32$ about $0.80$ of the prior probability is to the left of $0.05$ and
with $\beta=100$ about $0.99$ of the prior probability is to the left of
$0.05.$ We see that the misclassification rates for the small group $(c=1)$
stay about the same for $c_{RB}$ as $\beta$ increases while they deteriorate
markedly for $c_{MAP}$ as the MAP\ procedure basically ignores the small group.%

\begin{table}[tbp] \centering
\begin{tabular}
[c]{|r|r|r|}\hline
$\beta$ & $M_{0}(c_{MAP}\neq0)+M_{1}(c_{MAP}\neq1)$ & $M_{0}(c_{RB}%
\neq0)+M_{1}(c_{RB}\neq1)$\\\hline
$1$ & \multicolumn{1}{|c|}{$0.386+0.390=0.776$} &
\multicolumn{1}{|c|}{$0.386+0.390=0.776$}\\
$14$ & \multicolumn{1}{|c|}{$0.002+0.975=0.977$} &
\multicolumn{1}{|c|}{$0.285+0.380=0.665$}\\
$32$ & \multicolumn{1}{|c|}{$0.000+0.997=0.997$} &
\multicolumn{1}{|c|}{$0.292+0.349=0.641$}\\
$100$ & \multicolumn{1}{|c|}{$0.000+1.000=1.000$} &
\multicolumn{1}{|c|}{$0.300+0.324=0.624$}\\\hline
\end{tabular}
\caption{Conditional prior probabilities of misclassification for $c_{MAP}$ and $c_{RB} $ for various values of $\beta$ in Example 3 when $\alpha=1$, $\mu=1$, and $n$=10. }\label{TableKey}%
\end{table}%

We also investigated other choices for $n$ and $\mu.$ There is very little
change as $n$ increases. When $\mu$ moves towards $0$ the error rates go up
and go down as $\mu$ moves away from 0, as one would expect. It is the case,
however, that $c_{RB}$ always dominates $c_{MAP}.$ $\blacksquare$

\section{Estimation: continuous parameter space}

When $\psi$ has a continuous prior distribution the argument in Theorem 2 does
not work as $\Pi_{\Psi}(\{\delta(x)\}\,|\,x)=0.$ There are several possible
ways to proceed but one approach is to use a discretization of the problem
that uses Theorem 2.\ For this we will assume that the spaces involved are
locally Euclidean, mappings are sufficiently smooth and take the support
measures to be the analogs of Euclidean volume on the respective spaces. While
the argument provided applies quite generally, it is simplified here by taking
all spaces to be open subsets of Euclidean spaces and the support measures to
be Euclidean volume on these sets.

For each $\lambda>0$ suppose there is a discretization $\{B_{\lambda}%
(\psi):\psi\in\Psi(\Theta)\}$ of $\Psi(\Theta)$ into a countable number of
subsets with the following properties: $\psi\in B_{\lambda}(\psi),\Pi_{\Psi
}(B_{\lambda}(\psi))>0$ and $\sup_{\psi\in\Psi}$diam$(B_{\lambda}%
(\psi))\rightarrow0$ as $\lambda\rightarrow0.$ So, if $\psi^{\prime}\in
B_{\lambda}(\psi),$ then $B_{\lambda}(\psi^{\prime})=B_{\lambda}(\psi).$ For
example, the $B_{\lambda}(\psi)$ could be equal volume rectangles in $R^{k}.$
Further, we assume that $\Pi_{\Psi}(B_{\lambda}(\psi))/\nu_{\Psi}(B_{\lambda
}(\psi))\rightarrow\pi_{\Psi}(\psi)$ as $\lambda\rightarrow0$ for every
$\psi.$ This will hold whenever $\pi_{\Psi}$ is continuous everywhere and
$B_{\lambda}(\psi)$ converges nicely to $\{\psi\}$ as $\lambda\rightarrow0.$
Let $\psi_{\lambda}(\psi)$ denote a point in $B_{\lambda}(\psi)$ such that
$\psi_{\lambda}(\psi)=\psi_{\lambda}(\psi^{\prime})$ whenever $\psi
,\psi^{\prime}\in\in B_{\lambda}(\psi)$ and put $\Psi_{\lambda}=\{\psi
_{\lambda}(\psi):\psi\in\Psi(\Theta)\}.$ So $\Psi_{\lambda}$ is a discretized
version of $\Psi(\Theta).$ We will call this a \textit{regular discretization}
of $\Psi(\Theta).$ The discretized prior on $\Psi_{\lambda}$ is $\pi
_{\Psi,\lambda}(\psi_{\lambda}(\psi))=\Pi_{\Psi}(B_{\lambda}(\psi))$ and the
discretized posterior is $\pi_{\Psi,\lambda}(\psi_{\lambda}(\psi
)\,|\,x)=\Pi_{\Psi}(B_{\lambda}(\psi)\,|\,x).$

The loss function for the discretized problem is defined as in Theorem 2 by%
\begin{equation}
L_{RB,\lambda,\eta}(\theta,\psi_{\lambda}(\psi))=\frac{I_{\{\psi_{\lambda
}(\psi)\}}(\psi_{\lambda}(\Psi(\theta)))}{\max(\eta,\pi_{\Psi,\lambda}%
(\psi_{\lambda}(\Psi(\theta))))} \label{eq8}%
\end{equation}
and let $\delta_{\lambda,\eta}(x)$ denote a Bayes rule for this
problem.\smallskip

\noindent\textbf{Theorem 6}. Suppose that $\pi_{\Psi}$ is positive and
continuous and we have a regular discretization of $\Psi.$ Further suppose
that $\psi_{RB}(x)$ is the unique maximizer of $RB_{\Psi}(\psi\,|\,x)$ and for
any $\epsilon>0$\
\[
\sup_{\{\psi:||\psi-\psi_{RB}(x)||\geq\epsilon\}}RB_{\Psi}(\psi
\,|\,x)<RB_{\Psi}(\psi_{RB}(x)\,|\,x).
\]
Then, there exists $\eta(\lambda)\downarrow0$ as $\lambda\rightarrow0$ such
that a Bayes rule $\delta_{\lambda,\eta(\lambda)}(x),$ under the loss
$L_{RB,\lambda,\eta(\lambda)},$ converges to $\psi_{RB}(x)$ as $\lambda
\rightarrow0$ for all $x.$\smallskip

\noindent Theorem 6 says that $\psi_{RB}(x)$ is a limit of Bayes rules. So,
when $\Psi(\theta)=\theta$ we have the result that the MLE is a limit of Bayes
rules and more generally the MLE from an integrated likelihood is a limit of
Bayes rules. The regularity conditions stated in Theorem 6 hold in many common
statistical problems.

Now let $\hat{\psi}_{\lambda}(x)$ be the relative belief estimate from the
discretized problem, i.e., $\hat{\psi}_{\lambda}(x)$ maximizes $RB_{\Psi
}(B_{\lambda}(\psi)\,|\,x)$ as a function of $\psi\in\Psi_{\lambda}.$ The
following is immediate from the proof of Theorem 6, Theorem 5 and Corollary 4.
\smallskip

\noindent\textbf{Corollary 7}. $\hat{\psi}_{\lambda}$ is admissible and
Bayesian unbiased for the discretized problem and $\hat{\psi}_{\lambda
}(x)\rightarrow\psi_{RB}(x)$ as $\lambda\rightarrow0$ for every $x.$\smallskip

By similar arguments an analog of Theorem 6 for $\psi_{MAP}$ can be
established. Actually, in this case, a simpler development can be followed in
certain situations using the loss function $I_{B_{_{\lambda}}^{c}(\psi)}%
(\Psi(\theta))$. For this note that the posterior risk of $\delta$ in the
discretized problem is given by $1-\Pi_{\Psi}(B_{\lambda}(\delta
(x))\,|\,x)=1-\pi_{\Psi}(\delta^{\prime}(x)\,|\,x)\nu_{\Psi}(B_{\lambda
}(\delta(x)))$ for some $\delta^{\prime}(x)\in B_{\lambda}(\delta(x)).$ Now
suppose $B_{\lambda}(\psi)$ is a cube centered at $\psi$ of edge length
$\delta.$ Suppose further that for each $\epsilon>0$ there exists
$\lambda(\epsilon)>0$ such that, when $||\psi-\psi_{MAP}(x)||>\lambda
(\epsilon)$ then $\pi_{\Psi}(\psi\,|\,x)<\inf_{\psi^{\prime}\in B_{\lambda
(\epsilon)}(\psi_{MAP}(x))}\pi_{\Psi}(\psi^{\prime}\,|\,x).$ Since $\nu_{\Psi
}(B_{\lambda}(\psi))$ is constant we have that a Bayes rule $\delta
_{\lambda(\epsilon)}$ must then satisfy $||\delta_{\lambda(\epsilon)}%
(x)-\psi_{MAP}(x)||<\epsilon$. This proves that $\psi_{MAP}$ is a limit of
Bayes rules. By contrast, for the loss $I_{B_{_{\lambda}}^{c}(\psi)}%
(\Psi(\theta))/\Pi_{\Psi}(B_{\lambda}(\Psi(\theta))),$ the posterior risk of
$\delta$ is given by%
\[
\int_{\Psi}\{\Pi_{\Psi}(B_{\lambda}(\psi))\}^{-1}\Pi_{\Psi}(d\psi
\,|\,x)-\int_{B_{\lambda}(\delta(x))}\{\Pi_{\Psi}(B_{\lambda}(\psi))\}^{-1}%
\Pi_{\Psi}(d\psi\,|\,x)
\]
and the first term is generally unbounded unless $\Psi(\Theta)$ is compact.

Consider an important example.\smallskip

\noindent\textbf{Example 2.} \textit{Regression}

Suppose that $y=X\beta+e$ where $y\in R^{n},X\in R^{n\times k}$ is fixed of
rank $k,\beta\in R^{n\times k},$ and $e\sim N_{n}(0,\sigma^{2}I).$ We will
assume that $\sigma^{2}$ is known to simplify the discussion but this is not
necessary. Let $\pi$ be a prior density for $\beta.$ For every $\pi,$ having
observed $(X,y),$ then $\beta_{RB}(y)=b=(X^{\prime}X)^{-1}X^{\prime}y$ the MLE
of $\beta.$

It is interesting to contrast this result with what might be considered more
standard Bayesian estimates such as MAP or the posterior mean. For example,
suppose that $\beta\sim N_{k}(0,\tau^{2}I).$ Then the posterior distribution
of $\beta\ $is $N_{k}(\beta_{post}(y),\Sigma_{post})$ where
\[
\beta_{post}(y)=\Sigma_{post}(\sigma^{-2}X^{\prime}Xb),\text{ }\Sigma
_{post}=(\tau^{-2}I+\sigma^{-2}X^{\prime}X)^{-1}%
\]
and note $\beta_{MAP}(y)=\beta_{post}(y).$ Writing the spectral decomposition
of $X^{\prime}X$ as $X^{\prime}X=Q\Lambda Q^{\prime}$ we have that
\[
||\beta_{MAP}(y)||=||(I+(\sigma^{2}/\tau^{2})\Lambda^{-1})^{-1}Q^{\prime}b||.
\]
Since $||b||=||Q^{\prime}b||$ and $1/(1+\tau^{2}\lambda_{i}/\sigma^{2})<1$ for
each $i,$ this implies that $\beta_{MAP}(y)$ shrinks the MLE towards the prior
mean $0.$ When the columns of $X$ are orthonormal, then $\beta_{MAP}%
(y)=r(1+r)^{-1}b$ where $r=\tau^{2}/\sigma^{2}$ and so the shrinkage is
substantial unless $\tau^{2}$ is much larger than $\sigma^{2}.$ This shrinkage
is often cited as a positive attribute of these estimates. Consider, however,
the situation where the true value of $\beta$ is some distance from the mean.
In that case it seems wrong to move $\beta$ towards the prior mean and so it
isn't clear that shrinking the MLE is necessarily a good thing, particularly
as this requires giving up invariance.

Suppose it is required to estimate\ the mean response $\psi=\Psi
(\beta)=w^{\prime}\beta$ at $w$ for the predictors. The prior distribution of
$\psi$ is $N(0,\sigma_{\psi}^{2})=N(0,\tau^{2}w^{\prime}w)$ and the posterior
distribution is $N(\psi_{MAP}(y),\sigma_{\psi,post}^{2})=N(w^{\prime}%
\beta_{MAP}(y),w^{\prime}\Sigma_{post}(\beta)w).$ Note that
\[
\sigma_{\psi}^{2}-\sigma_{\psi,post}^{2}=w^{\prime}(\tau^{2}I-\Sigma
_{post})w=\tau^{2}w^{\prime}Q^{\prime}(I-(I+(\tau^{2}/\sigma^{2})\Lambda
)^{-1})Qw>0
\]
since $1/(1+\tau^{2}\lambda_{i}/\sigma^{2})<1$ for each $i.$ Therefore,
maximizing the ratio of the posterior to prior densities leads to%
\begin{equation}
\psi_{RB}(y)=(1-\sigma_{\psi,post}^{2}/\sigma_{\psi}^{2})^{-1}\psi_{MAP}(y).
\label{est1}%
\end{equation}
Then $\sigma_{\psi}^{2}>\sigma_{\psi,post}^{2}$ implies $|\psi_{RB}%
(y)|>|\psi_{MAP}(y)|.$ Note that when $\sigma_{\psi,post}^{2}$ is much smaller
than $\sigma_{\psi}^{2},$ in other words the posterior is much more
concentrated than the prior, then $\psi_{RB}(y)$ and $\psi_{MAP}(y)$ are very
similar. In general $\psi_{RB}(y)$ is not equal to $w^{\prime}b,$ the plug-in
MLE of $\psi,$ although it is the MLE from the integrated likelihood,
$\psi_{RB}(y)\rightarrow w^{\prime}b$ as $\tau^{2}\rightarrow\infty$ and when
$X$ has orthonormal columns $\psi_{RB}(y)=w^{\prime}b.$

Suppose it is required to predict a response $z$ at the predictor value $w\in
R^{k}.$ When $\beta\sim N_{k}(0,\tau^{2}I)$ the prior distribution of $z$ is
$z\sim N(0,\sigma^{2}+\tau^{2}w^{\prime}w)=N(0,\sigma_{z}^{2})$ and the
posterior distribution is $N(\mu_{post}(z),\sigma_{post}^{2}(z))$ where
\[
\mu_{post}(z)=w^{\prime}\beta_{post}(y),\text{ }\sigma_{post}^{2}%
(z)=\sigma^{2}+w^{\prime}\Sigma_{post}w.
\]
To obtain $z_{RB}(y)$ it is necessary to maximize the ratio of the posterior
to the prior densities of $z$ and this leads to%
\begin{equation}
z_{RB}(y)=(1-\sigma_{post}^{2}(z)/\sigma_{prior}^{2}(z))^{-1}\mu_{post}(z).
\label{pred1}%
\end{equation}
Note that $\sigma_{z}^{2}-\sigma_{post}^{2}(z)=\sigma_{z}^{2}(w^{\prime}%
\beta)-\sigma_{post}^{2}(w^{\prime}\beta)>0$ and so $|z_{RB}(y)|>|\mu
_{post}(z)|$ and $z_{RB}$ is further from the prior mean than $z_{MAP}%
(y)=\mu_{post}(z).$ Also, we see that, when $\sigma_{post}^{2}(z)$ is small
then $z_{RB}(y)$ and $z_{MAP}(y)$ are very similar. Finally, comparing
(\ref{est1}) and (\ref{pred1}) we have that
\[
z_{RB}(y)=(\sigma_{prior}^{2}(z)/\sigma_{post}^{2}(\psi))w^{\prime}\psi
_{RB}(y)=(1+\sigma^{2}/\tau^{2})\psi_{RB}(y)
\]
and so $\psi_{RB}(y)$ at $w$ is more dispersed than the $z_{RB}(y)$ estimate
of the mean at $w$ and this makes good sense as we have to take into account
the additional variation due to prediction. By contrast $w_{MAP}(y)=\psi
_{MAP}(y).$ $\blacksquare$

\section{Credible regions and hypothesis assessment}

First recall that a $\gamma$-relative belief credible region for $\psi
=\Psi(\theta)$ is given by $C_{\Psi,\gamma}(x)=\{\psi:RB_{\Psi}(\psi
\,|\,x)\geq c_{\gamma}(x)\}$ where $c_{\gamma}(x)=\sup\{c:\Pi_{\Psi}(RB_{\Psi
}(\psi\,|\,x)\geq c\,|\,x)\geq\gamma\}.$ There is some arbitrariness in the
choice of the greater than or equal sign to define the credible region as it
also could have been defined as $C_{\Psi,\gamma}(x)=\{\psi:RB_{\Psi}%
(\psi\,|\,x)>c_{\gamma}(x)\}$ where $c_{\gamma}(x)=\inf\{c:\Pi_{\Psi}%
(RB_{\Psi}(\psi\,|\,x)\leq c\,|\,x)\leq1-\gamma\}.$ In this latter case
$c_{\gamma}(x)$ is the $(1-\gamma)$-th quantile of the posterior distribution
of the relative belief ratio$.$ This definition has some advantages as using
this implies that the plausible region satisfies $Pl_{\Psi}(x)=C_{\Psi,\gamma
}(x)$ where $\gamma=\Pi_{\Psi}(Pl_{\Psi}(x)\,|\,x).$ Also, the strength of the
evidence concerning the hypothesis $H_{0}:\Psi(\theta)=\psi_{0}$ satisfies
$Str_{\Psi}(\psi_{0}\,|\,x)=1-\Pi_{\Psi}(C_{\Psi,\gamma}(x)\,|\,x)$ where
$\gamma=1-Str_{\Psi}(\psi_{0}\,|\,x).$ The point here is that there is a close
relationship between relative belief credible regions and the plausible region
and the strength calculation. As such, any decision-theoretic interpretation
for relative belief credible regions also applies to the plausible region and
the strength of the evidence. Throughout this section we will, however, retain
the definition for $C_{\Psi,\gamma}(x)$ provided in Section 2.3.

Now consider the lowest posterior loss $\gamma$-credible regions that arise
from the prior-based loss functions considered here.\smallskip

\noindent\textbf{Theorem 8}. Suppose that $\pi_{\Psi}(\psi)>0$ for every
$\psi\in\Psi(\Theta)$ where $\Psi(\Theta)$ is finite with $\nu_{\Psi}$ equal
to counting measure. Then $C_{\Psi,\gamma}(x)$ is a $\gamma$-lowest posterior
loss credible region for the loss function $L_{RB}.$

\noindent Proof: From (\ref{eq3}) and (\ref{eq5a}) the $\gamma$-lowest
posterior loss credible region is%
\[
D_{\gamma}(x)=\left\{  \psi:RB_{\Psi}(\psi\,|\,x)\geq\int_{\Psi}RB_{\Psi
}(\zeta\,|\,x)\,\nu_{\Psi}(d\zeta)-d_{\gamma}(x)\right\}
\]
and $d_{\gamma}(x)=\sup\{d:\Pi_{\Psi}(r(\psi\,|\,x)\leq d\,|\,x)\geq\gamma\}.$
As $\int_{\Psi}RB_{\Psi}(z\,|\,x)\,\nu_{\Psi}(dz)$ is independent of $\psi$ it
is clearly equivalent to define this region via $C_{\Psi,\gamma}(x)=\left\{
\psi:RB_{\Psi}(\psi\,|\,x)\geq c_{\gamma}(x)\right\}  ,$ namely, $D_{\gamma
}(x)=C_{\Psi,\gamma}(x).$ $\blacksquare$\smallskip

Now consider the case where $\Psi$ is countable and we use loss function
$L_{RB,\eta}.$ Following the proof of Theorem 8 we see that a $\gamma$-lowest
posterior loss region takes the form%
\[
D_{\eta,\gamma}(x)=\left\{  \psi:\pi_{\Psi}(\psi\,|\,x)/\max(\eta,\pi_{\Psi
}(\psi))\geq d_{\eta,\gamma}(x)\right\}
\]
where $d_{\eta,\gamma}(x)=\sup\{d:\Pi_{\Psi}(\pi_{\Psi}(\psi\,|\,x)/\max
(\eta,\pi_{\Psi}(\psi))\,|\,x)\geq d)\geq\gamma\}.$\smallskip

\noindent\textbf{Theorem 9}. Suppose that $\pi_{\Psi}(\psi)>0$ for every
$\psi\in\Psi,$ that $\Psi$ is countable with $\nu_{\Psi}$ equal to counting
measure. For the loss function $L_{RB,\eta}$ then $C_{\Psi,\gamma}%
(x)\subset\lim\inf_{\eta\rightarrow0}D_{\eta,\gamma}(x)$ whenever $\gamma$ is
such that $\Pi_{\Psi}(C_{\Psi,\gamma}(x)\,|\,x)=\gamma$ and $\lim\sup
_{\eta\rightarrow0}D_{\eta,\gamma}(x)\subset C_{\Psi,\gamma^{\prime}}(x)$
whenever $\gamma^{\prime}>\gamma$ and $\Pi_{\Psi}(C_{\Psi,\gamma^{\prime}%
}(x)\,|\,x)=\gamma^{\prime}.$\smallskip

\noindent While Theorem 9 does not establish the exact convergence $\lim
_{\eta\rightarrow0}D_{\eta,\gamma}(x)=C_{\gamma}(x)$ it is likely, however,
that this does hold under quite general circumstances due to the discreteness.
Theorem 9 does show that limit points of the class of sets $D_{\eta,\gamma
}(x)$ always contain $C_{\Psi,\gamma}(x)$ and their posterior probability
content differs from $\gamma$ by at most $\gamma^{\prime}-\gamma$ where
$\gamma^{\prime}>\gamma$ is the next largest value for which we have exact content.

Now consider the continuous case with a regular discretization. For $S^{\ast
}\subset\Psi_{\lambda}=\{\psi_{\lambda}(\psi):\psi_{\lambda}(\psi)\in
B_{\lambda}(\psi)\},$ namely, $S^{\ast}$ is a subset of a discretized version
of $\Psi(\Theta),$ define the \textit{undiscretized} version of $S^{\ast}$ to
be $S=\cup_{\psi\in S^{\ast}}B_{\lambda}(\psi).$ Now let $C_{\Psi
,\lambda,\gamma}^{\ast}(x)$ be the $\gamma$-relative belief region for the
discretized problem and let $C_{\Psi,\lambda,\gamma}(x)$ be its undiscretized
version. Note that in a continuous context we will consider two sets as equal
if they differ only by a set of measure 0 with respect to $\Pi_{\Psi}.$ The
following result says that a $\gamma$-relative belief credible region for the
discretized problem, after undiscretizing, converges to the $\gamma$-relative
belief region for the original problem.\smallskip

\noindent\textbf{Theorem 10}. Suppose that $\pi_{\Psi}$ is positive and
continuous, there is a regular discretization of $\Psi(\Theta)$ and $RB_{\Psi
}(\psi\,|\,x)$ has a continuous posterior distribution. Then $\lim
_{\lambda\rightarrow0}C_{\Psi,\lambda,\gamma}(x)=C_{\Psi,\gamma}%
(x).$\smallskip

\noindent While Theorem 10 has interest in its own right, it can be also used
to prove that relative belief regions are limits of lowest posterior loss regions.

Let $D_{\eta,\lambda,\gamma}^{\ast}(x)$ be the $\gamma$-lowest posterior loss
region obtained for the discretized problem using loss function (\ref{eq8})
and $D_{\eta,\lambda,\gamma}(x)$ be the undiscretized version. \smallskip

\noindent\textbf{Theorem 11}. Suppose that $\pi_{\Psi}$ is positive and
continuous, we have a regular discretization of $\Psi$ and $RB_{\Psi}%
(\psi\,|\,x)$ has a continuous posterior distribution. Then $C_{\Psi,\gamma
}(x)=\underset{\lambda\rightarrow0}{\lim}\underset{\eta\rightarrow0}{\lim\inf
}D_{\Psi,\gamma}(x)=\underset{\lambda\rightarrow0}{\lim}\underset
{\eta\rightarrow0}{\lim\sup}D_{\Psi,\gamma}(x).$

In Evans, Guttman, and Swartz (2006) and Evans and Shakhatreh (2008)
additional properties of relative belief regions are developed. For example,
it is proved that a $\gamma$-relative belief region $C_{\Psi,\gamma}(x)$ for
$\psi$ satisfying $\Pi_{\Psi}(C_{\Psi,\gamma}(x)\,|\,x)=\gamma$ minimizes
$\Pi_{\Psi}(B)$ among all (measurable) subsets of $\Psi$ satisfying $\Pi
_{\Psi}(B\,|\,x)\geq\gamma.$ So a $\gamma$-relative belief region is smallest
among all $\gamma$-credible regions for $\psi$ where size is measured using
the prior measure. This property has several consequences. For example, the
prior probability that a region $B(x)\subset\Psi(\Theta)$ contains a false
value from the prior is given by $\int_{\Theta}\int_{\Psi}F_{\theta}(\psi\in
B(x))\,\Pi_{\Psi}(d\psi)\,\Pi(d\theta)$ where a false value is a value of
$\psi\sim\Pi_{\Psi}$ generated independently of $(\theta,x)\sim\Pi_{\Psi
}\times F_{\theta}.$ It can be proved that a $\gamma$-relative belief region
minimizes this probability among all $\gamma$-credible regions for $\psi$ and
is always unbiased in the sense that the probability of covering a false value
is bounded above by $\gamma.$ Furthermore, a $\gamma$-relative belief region
maximizes the relative belief ratio $\Pi_{\Psi}(B\,|\,x)/\Pi_{\Psi}(B)$ and
the Bayes factor $\Pi_{\Psi}(B\,|\,x)\Pi_{\Psi}(B^{c})/\Pi_{\Psi}%
(B^{c}\,|\,x)\Pi_{\Psi}(B)$ among all regions $B\subset\Psi$ with $\Pi_{\Psi
}(B)=\Pi_{\Psi}(C_{\Psi,\gamma}(x)\,|\,x).$

While the results in this section have been concerned with obtaining credible
regions for parameters, similar results can be proved for the construction of
prediction regions.

\section{Conclusions}

Relative belief inferences are closely related to likelihood inferences. This
together with their invariance and optimality properties make these prime
candidates as appropriate inferences in Bayesian contexts.\ This paper has
shown that relative belief inferences arise naturally in a decision-theoretic
formulation using loss functions based on the prior.

\section*{Appendix}

\noindent\textbf{Proof of Theorem} \textbf{2 and Corollary 3}: We have that
\begin{align}
r_{\eta}(\delta\,|\,x)  &  =\int_{\Psi}L_{RB,\eta}(\theta,\delta(x))\pi_{\Psi
}(\psi\,|\,x)\,\nu_{\Psi}(d\psi)\nonumber\\
&  =\int_{\Psi}\frac{\pi_{\Psi}(\psi\,|\,x)}{\max(\eta,\pi_{\Psi}(\psi))}%
\,\nu_{\Psi}(d\psi)-\frac{\pi_{\Psi}(\delta(x)\,|\,x)}{\max(\eta,\pi_{\Psi
}(\delta(x)))}. \label{eq6a}%
\end{align}
The first term in (\ref{eq6a}) is bounded above by $1/\eta$ and does not
depend on $\delta(x)$ so the value of a Bayes rule at $x$ is obtained by
finding $\delta(x)$ that maximizes the second term. Note that
\begin{equation}
\frac{\pi_{\Psi}(\delta(x)\,|\,x)}{\max(\eta,\pi_{\Psi}(\delta(x)))}=\left\{
\begin{array}
[c]{cl}%
\frac{\pi_{\Psi}(\delta(x)\,|\,x)}{\eta} & \text{if }\eta>\pi_{\Psi}%
(\delta(x)),\\
RB_{\Psi}(\delta(x)\,|\,x) & \text{if }\eta\leq\pi_{\Psi}(\delta(x)).
\end{array}
\right.  \label{eq6b}%
\end{equation}
There are at most finitely many values of $\psi$ satisfying $\eta\leq\pi
_{\Psi}(\psi)$ and so $RB_{\Psi}(\psi\,|\,x)$ assumes a maximum on this set,
say at $\psi_{\eta}(x),$ and $\psi_{\eta}(x)=\psi_{RB}(x)$ when $\eta\leq
\pi_{\Psi}(\psi_{RB}(x)).$ If $\eta>\pi_{\Psi}(\delta(x)),$ then $\pi_{\Psi
}(\delta(x)\,|\,x)/\eta<RB_{\Psi}(\delta(x)\,|\,x)\leq RB_{\Psi}(\psi
_{RB}(x)\,|\,x).$ This proves that, for all $\eta\leq\eta_{0}=\pi_{\Psi}%
(\psi_{RB}(x))>0$ the maximizer of (\ref{eq6b}) is given by $\delta
(x)=\psi_{RB}(x)$ and the results are established.\medskip

\noindent\textbf{Proof of Theorem 5}: The prior risk of $\delta$ is given by
\begin{align*}
r(\delta)  &  =\int_{\Theta}\int_{\mathcal{X}}L(\theta,\delta
(x))\,F(dx\,|\,\theta)\,\Pi(d\theta)\\
&  =\int_{\Theta}\int_{\mathcal{X}}[h(\Psi(\theta))-I_{\{\delta(x)\}}%
(\Psi(\theta))h(\Psi(\theta))]\,F(dx\,|\,\theta)\,\Pi(d\theta)\\
&  =\int_{\Theta}h(\Psi(\theta))\,\Pi(d\theta)-\int_{\mathcal{X}}\int_{\Theta
}I_{\{\delta(x)\}}(\Psi(\theta))h(\Psi(\theta))\,\Pi(d\theta\,|\,x)\,M(dx)\\
&  =\int_{\Theta}h(\Psi(\theta))\,\Pi(d\theta)-\int_{\mathcal{X}}%
h(\delta(x))\pi_{\Psi}(\delta(x)\,|\,x)\,M(dx)
\end{align*}
and%
\begin{align*}
&  \int_{\Theta}\int_{\Theta}\int_{\mathcal{X}}L(\theta^{\prime}%
,\delta(x))\,F(dx\,|\,\theta)\,\Pi(d\theta)\,\Pi(d\theta^{\prime})\\
&  =\int_{\Theta}\int_{\Theta}\int_{\mathcal{X}}[h(\Psi(\theta^{\prime
}))-I_{\{\delta(x)\}}(\Psi(\theta^{\prime}))h(\Psi(\theta^{\prime
}))]\,F(dx\,|\,\theta)\,\Pi(d\theta)\,\Pi(d\theta^{\prime})\\
&  =\int_{\Theta}h(\Psi(\theta))\,\Pi(d\theta)-\int_{\mathcal{X}}%
h(\delta(x))\pi_{\Psi}(\delta(x))\,M(dx).
\end{align*}
Therefore, $\delta$ is Bayesian unbiased if and only if
\begin{equation}
\int_{\mathcal{X}}h(\delta(x))[\pi_{\Psi}(\delta(x)\,|\,x)-\pi_{\Psi}%
(\delta(x))]\,M(dx)\geq0. \label{eq6d}%
\end{equation}
This inequality holds when $\delta(x)=\psi_{RB}(x)$ because $\pi_{\Psi}%
(\cdot\,|\,x)/\pi_{\Psi}(\cdot)$ is the density of $\Pi_{\Psi}(\cdot\,|\,x)$
with respect to $\Pi_{\Psi}$ and which implies that the maximum of this
density is greater than or equal to 1. \medskip

\noindent\textbf{Proof of Theorem 6 and Corollary 7}: Just as in Theorem 2, a
Bayes rule $\delta_{\lambda,\eta}(x)$ maximizes $\pi_{\Psi,\lambda}%
(\delta(x)\,|\,x)/\max(\eta,\pi_{\Psi,\lambda}(\delta(x)))$ for $\delta
(x)\in\Psi_{\lambda}.$ Furthermore, as in Theorem 2, such a rule exists. Now
define $\eta(\lambda)$ so that $0<\eta(\lambda)<\Pi_{\Psi}(B_{\lambda}%
(\psi_{RB}(x)))$ and note that $\eta(\lambda)\rightarrow0$ as $\lambda
\rightarrow0.$ We have that, as $\lambda\rightarrow0,$
\begin{equation}
\frac{\pi_{\Psi,\lambda}(\psi_{\lambda}(\psi_{RB}(x))\,|\,x)}{\max
(\eta(\lambda),\pi_{\Psi,\lambda}(\psi_{\lambda}(\psi_{RB}(x)))}=\frac
{\pi_{\Psi,\lambda}(\psi_{\lambda}(\psi_{RB}(x))\,|\,x)}{\pi_{\Psi,\lambda
}(\psi_{\lambda}(\psi_{RB}))}\rightarrow RB_{\Psi}(\psi_{RB}(x)\,|\,x).
\label{eq7a}%
\end{equation}
Let $\epsilon>0.$ Let $\lambda_{0}$ be such that $\sup_{\psi\in\Psi}%
$diam$(B_{\lambda}(\psi))<\epsilon/2$ for all $\lambda<\lambda_{0}.$ Then for
$\lambda<\lambda_{0},$ and any $\delta(x)$ satisfying $||\delta(x)-\psi
_{RB}(x)||\geq\epsilon,$ we have
\begin{align}
&  \frac{\pi_{\Psi,\lambda}(\psi_{\lambda}(\delta(x))\,|\,x)}{\pi
_{\Psi,\lambda}(\psi_{\lambda}(\delta(x)))}=\frac{\int_{B_{\lambda}%
(\psi_{\lambda}(\delta(x)))}\pi_{\Psi}(\psi\,|\,x)\,\nu_{\Psi}(d\psi)}%
{\int_{B_{\lambda}(\psi_{\lambda}(\delta(x)))}\pi_{\Psi}(\psi)\,\nu_{\Psi
}(d\psi)}\nonumber\\
&  =\frac{\int_{B_{\lambda}(\psi_{\lambda}(\delta(x)))}RB_{\Psi}%
(\psi\,|\,x)\pi_{\Psi}(\psi)\,\nu_{\Psi}(d\psi)}{\int_{B_{\lambda}%
(\psi_{\lambda}(\delta(x)))}\pi_{\Psi}(\psi)\,\nu_{\Psi}(d\psi)}\nonumber\\
&  \leq\sup_{\{\psi:||\psi-\psi_{RB}(x)||>\epsilon/2\}}RB_{\Psi}%
(\psi\,|\,x)<RB_{\Psi}(\psi_{RB}(x)\,|\,x). \label{eq7b}%
\end{align}
\ By (\ref{eq7a}) and (\ref{eq7b}) there exists $\lambda_{1}<\lambda_{0}$ such
that, for all\ $\lambda<\lambda_{1},$%
\begin{equation}
\frac{\pi_{\Psi,\lambda}(\psi_{\lambda}(\psi_{RB}(x))\,|\,x)}{\pi
_{\Psi,\lambda}(\psi_{\lambda}(\psi_{RB}(x)))}>\sup_{\{\psi:||\psi-\psi
_{RB}(x)||>\epsilon/2\}}RB_{\Psi}(\psi\,|\,x). \label{eq7c}%
\end{equation}
Therefore, when $\lambda<\lambda_{1},$ a Bayes rule $\delta_{\lambda
,\eta(\lambda)}(x)$ satisfies
\begin{align}
&  \frac{\pi_{\Psi,\lambda}(\delta_{\lambda,\eta(\lambda)}(x)\,|\,x)}%
{\pi_{\Psi,\lambda}(\delta_{\lambda,\eta(\lambda)}(x))}\geq\frac{\pi
_{\Psi,\lambda}(\delta_{\lambda,\eta(\lambda)}(x)\,|\,x)}{\max(\eta
(\lambda),\pi_{\Psi,\lambda}(\delta_{\lambda,\eta(\lambda)}(x)))}\nonumber\\
&  \geq\frac{\pi_{\Psi,\lambda}(\psi_{\lambda}(\psi_{RB}(x))\,|\,x)}{\max
(\eta(\lambda),\pi_{\Psi,\lambda}(\psi_{\lambda}(\psi_{RB}(x)))}=\frac
{\pi_{\Psi,\lambda}(\psi_{\lambda}(\psi_{RB}(x))\,|\,x)}{\pi_{\Psi,\lambda
}(\psi_{\lambda}(\psi_{RB}(x)))}. \label{eq7d}%
\end{align}
By (\ref{eq7b}), (\ref{eq7c}) and (\ref{eq7d}) this implies that
$||\delta_{\lambda,\eta(\lambda)}(x)-\psi_{RB}(x)||<\epsilon/2$ and the
convergence is established.

Now $\pi_{\Psi,\lambda}(\hat{\psi}_{\lambda}(x)\,|\,x)/\pi_{\Psi,\lambda}%
(\hat{\psi}_{\lambda}(x))\geq\pi_{\Psi,\lambda}(\delta_{\lambda,\eta(\lambda
)}(x)\,|\,x)/\pi_{\Psi,\lambda}(\delta_{\lambda,\eta(\lambda)}(x))$ and so by
(\ref{eq7b}), (\ref{eq7c}) and (\ref{eq7d}) this implies that $||\hat{\psi
}_{\lambda}(x)-\psi_{RB}(x)||<\epsilon$ and the convergence of $\hat{\psi
}_{\lambda}(x)\ $to $\psi_{RB}(x)$ is established.\medskip

\noindent\textbf{Proof of Theorem 9}: For $c>0$ let $S_{c}(x)=\{\psi:RB_{\Psi
}(\psi\,|\,x)\geq c\}$ and $S_{\eta,c}(x)=\{\psi:\pi_{\Psi}(\psi
\,|\,x)/\max(\eta,\pi_{\Psi}(\psi))\geq c\}.$ Note that $S_{\eta,c}(x)\uparrow
S_{c}(x)$ as $\eta\rightarrow0.$

Suppose $c$ is such that $\Pi_{\Psi}(S_{c}(x)\,|\,x)\leq\gamma.$ Then
$\Pi_{\Psi}(S_{\eta,c}(x)\,|\,x)\leq\gamma$ for all $\eta$ and so $S_{\eta
,c}(x)\subset D_{\eta,\gamma}(x).$ This implies that $S_{c}(x)\subset\lim
\inf_{\eta\rightarrow0}D_{\eta,\gamma}(x)$ and since $\Pi_{\Psi}%
(C_{\Psi,\gamma}(x)\,|\,x)=\gamma$ this implies that $C_{\Psi,\gamma
}(x)\subset\lim\inf_{\eta\rightarrow0}D_{\eta,\gamma}(x).$

Now suppose $c$ is such that $\Pi_{\Psi}(S_{c}(x)\,|\,x)>\gamma.$ Then there
exists $\eta_{0}$ such that for all $\eta<\eta_{0}$ we have $\Pi_{\Psi
}(S_{\eta,c}(x)\,|\,x)>\gamma.$ Since $D_{\eta,\gamma}(x)\subset S_{\eta
,c}(x)$ when $\eta<\eta_{0},$ then $\lim\sup_{\eta\rightarrow0}D_{\eta,\gamma
}(x)\subset S_{c}(x).$ Then choosing $c=c_{\gamma^{\prime}}(x)$ for
$\gamma^{\prime}>\gamma$ implies that $\lim\sup_{\eta\rightarrow0}%
D_{\eta,\gamma}(x)\subset C_{\Psi,\gamma^{\prime}}(x).$\medskip

\noindent\textbf{Proof of Theorem 10}: Let $S_{c}(x)=\{\psi:RB_{\Psi}%
(\psi\,|\,x)\geq c\}$ and $S_{\lambda,c}(x)=\{\psi:\Pi_{\Psi}(B_{\lambda}%
(\psi)\,|\,x)/\Pi_{\Psi}(B_{\lambda}(\psi))\geq c\}.$ Recall that
\[
\lim_{\lambda\rightarrow0}\Pi_{\Psi}(B_{\lambda}(\psi)\,|\,x)/\Pi_{\Psi
}(B_{\lambda}(\psi))=\lim_{\lambda\rightarrow0}RB_{\Psi}(B_{\lambda}%
(\psi)\,|\,x)=RB_{\Psi}(\psi\,|\,x)
\]
for every $\psi.$ If $RB_{\Psi}(\psi\,|\,x)>c,$ there exists $\lambda_{0}$
such that for all $\lambda<\lambda_{0},$ then $\Pi_{\Psi}(RB_{\Psi}%
(B_{\lambda}(\psi)\,|\,x)>c$ and this implies that $\psi\in\lim\inf
_{\lambda\rightarrow0}S_{\lambda,c}(x).$ Now $\Pi_{\Psi}(RB_{\Psi}%
(\psi\,|\,x)=c)=0$ and so $S_{c}(x)\subset\lim\inf_{\lambda\rightarrow
0}S_{\lambda,c}(x)$ (after possibly deleting a set of $\Pi_{\Psi}$-measure 0
from $S_{c}(x)).$ If $\psi\in\lim\sup_{\lambda\rightarrow0}S_{\lambda,c}(x),$
then $RB_{\Psi}(B_{\lambda}(\psi)\,|\,x)\geq c$ for infinitely many
$\lambda\rightarrow0,$ which implies that $RB_{\Psi}(\psi\,|\,x)\geq c,$ and
therefore $\psi\in$ $S_{c}(x).$ This proves $S_{c}(x)=\lim_{\lambda
\rightarrow0}S_{\lambda,c}(x)$ (up to a set of $\Pi_{\Psi}$-measure 0) so that
$\lim_{\lambda\rightarrow0}\Pi_{\Psi}(S_{\lambda,c}(x)\Delta S_{c}%
(x)\,|\,x)=0$ for any $c.$

Let $c_{\lambda,\gamma}(x)=\sup\{c\geq0:\Pi_{\Psi}(S_{\lambda,c}%
(x)\,|\,x)\geq\gamma\}$ so $S_{c_{\gamma}(x)}(x)=C_{\Psi,\gamma}(x),$
$S_{\lambda,c_{\lambda,\gamma}(x)}(x)=C_{\Psi,\lambda,\gamma}(x)\ $and
\begin{align}
&  \Pi_{\Psi}(C_{\Psi,\gamma}(x)\Delta C_{\Psi,\lambda,\gamma}(x)\,|\,x)=\Pi
_{\Psi}(S_{c_{\gamma}(x)}(x)\Delta S_{\lambda,c_{\lambda,\gamma}%
(x)}(x)\,|\,x)\nonumber\\
&  \leq\Pi_{\Psi}(S_{c_{\gamma}(x)}(x)\Delta S_{\lambda,c_{\gamma}%
(x)}(x)\,|\,x)+\Pi_{\Psi}(S_{\lambda,c_{\lambda,\gamma}(x)}(x)\Delta
S_{\lambda,c_{\gamma}(x)}(x)\,|\,x). \label{ineq}%
\end{align}
Since $S_{c_{\gamma}(x)}(x)=\lim_{\lambda\rightarrow0}S_{\lambda,c_{\gamma
}(x)}(x)$ then $\Pi_{\Psi}(S_{c_{\gamma}(x)}(x)\Delta S_{\lambda,c_{\gamma
}(x)}(x)\,|\,x)\rightarrow0$ and $\Pi_{\Psi}(S_{\lambda,c_{\gamma}%
(x)}(x)\,|\,x)\rightarrow\Pi_{\Psi}(S_{c_{\gamma}(x)}(x)\,|\,x)=\gamma$ as
$\lambda\rightarrow0.$ Now consider the second term in (\ref{ineq}). Since
$RB_{\Psi}(\psi\,|\,x)$ has a continuous posterior distribution, $\Pi_{\Psi
}(RB_{\Psi}(\psi\,|\,x)\geq c\,|\,x)$ is continuous in $c.$ Let $\epsilon>0$
and note that for all $\lambda$ small enough, $\Pi_{\Psi}(S_{\lambda
,c_{\gamma-\epsilon}(x)}(x)\,|\,x)<\gamma$ and $\Pi_{\Psi}(S_{\lambda
,c_{\gamma+\epsilon}(x)}(x)\,|\,x)>\gamma$ which implies that $c_{\gamma
+\epsilon}(x)\leq c_{\lambda,\gamma}(x)\leq c_{\gamma-\epsilon}(x)$ and
therefore $S_{\lambda,c_{\gamma+\epsilon}(x)}(x)\subset S_{\lambda
,c_{\lambda,\gamma}(x)}\subset S_{\lambda,c_{\gamma-\epsilon}(x)}(x).$ As
$S_{\lambda,c_{\lambda,\gamma}(x)}(x)\subset S_{\lambda,c_{\gamma}(x)}(x)$ or
$S_{\lambda,c_{\lambda,\gamma}(x)}(x)\supset S_{\lambda,c_{\gamma}(x)}(x)$
then
\[
\Pi_{\Psi}(S_{\lambda,c_{\lambda,\gamma}(x)}(x)\Delta S_{\lambda,c_{\gamma
}(x)}(x)\,|\,x)=|\Pi_{\Psi}(S_{\lambda,c_{\lambda,\gamma}(x)}(x)\,|\,x)-\Pi
_{\Psi}(S_{\lambda,c_{\gamma}(x)}(x)\,|\,x)|.
\]
For all $\lambda$ small, then $|\Pi_{\Psi}(S_{\lambda,c_{\lambda,\gamma}%
(x)}(x)\,|\,x)-\Pi_{\Psi}(S_{\lambda,c_{\gamma}(x)}(x)\,|\,x)|$ is bounded
above by%
\begin{align*}
&  \max\{|\Pi_{\Psi}(S_{\lambda,c_{\gamma+\epsilon}(x)}(x)\,|\,x)-\Pi_{\Psi
}(S_{\lambda,c_{\gamma}(x)}(x)\,|\,x)|,\\
&  \qquad|\Pi_{\Psi}(S_{\lambda,c_{\gamma-\epsilon}(x)}(x)\,|\,x)-\Pi_{\Psi
}(S_{\lambda,c_{\gamma}(x)}(x)\,|\,x)|\}
\end{align*}
and this upper bound converges to $\epsilon$ as $\lambda\rightarrow0.$ Since
$\epsilon$ is arbitrary this implies that the second term in (\ref{ineq}) goes
to 0 as $\lambda\rightarrow0$ and this proves the result.\medskip

\noindent\textbf{Proof of Theorem 11}: Suppose, without loss of generality
that $0<\gamma<1.$ Let $\epsilon>0$ and $\delta>0$ satisfy $\gamma+\delta
\leq1.$ Put $\gamma^{\prime}(\lambda,\gamma)=\Pi_{\Psi}(C_{\Psi,\lambda
,\gamma}(x)\,|\,x),\gamma^{\prime\prime}(\lambda,\gamma)$\newline$=\Pi_{\Psi
}(C_{\Psi,\lambda,\gamma+\delta}(x)\,|\,x)$ and note that $\gamma^{\prime
}(\lambda,\gamma)\geq\gamma,\gamma^{\prime\prime}(\lambda,\gamma)\geq
\gamma+\delta.$ By Theorem 10 $C_{\Psi,\lambda,\gamma}(x)\rightarrow
C_{\Psi,\gamma}(x)$ and $C_{\Psi,\lambda,\gamma+\delta}(x)\rightarrow
C_{\Psi,\gamma+\delta}(x)$ as $\lambda\rightarrow0$ so $\gamma^{\prime
}(\lambda,\gamma)\rightarrow\gamma$ and $\gamma^{\prime\prime}(\lambda
,\gamma)\rightarrow\gamma+\delta$ as $\lambda\rightarrow0.$ This implies that
there is a $\lambda_{0}(\delta)$ such that for all $\lambda<\lambda_{0}%
(\delta)$ then $\gamma^{\prime}(\lambda,\gamma)<\gamma^{\prime\prime}%
(\lambda,\gamma).$ Therefore, by Theorem 9, we have that for all
$\lambda<\lambda_{0}(\delta)$%
\begin{equation}
C_{\Psi,\lambda,\gamma}(x)\subset\underset{\eta\rightarrow0}{\lim\inf
\,}D_{\eta,\lambda,\gamma^{\prime}(\lambda,\gamma)}(x)\subset\underset
{\eta\rightarrow0}{\lim\sup\,}D_{\eta,\lambda,\gamma^{\prime}(\lambda,\gamma
)}(x)\subset C_{\Psi,\lambda,\gamma+\delta}(x). \label{set1}%
\end{equation}
From (\ref{set1}) and Theorem 10 we have that
\begin{align*}
C_{\Psi,\gamma}(x)  &  \subset\underset{\lambda\rightarrow0}{\lim\inf
\,}\underset{\eta\rightarrow0}{\lim\inf\,}D_{\eta,\lambda,\gamma^{\prime
}(\lambda,\gamma)}(x)\\
&  \subset\underset{\lambda\rightarrow0}{\lim\sup}\underset{\eta\rightarrow
0}{\,\lim\sup\,}D_{\eta,\lambda,\gamma^{\prime}(\lambda,\gamma)}(x)\subset
C_{\Psi,\gamma+\delta}(x).
\end{align*}
Since $\lim_{\delta\rightarrow0}C_{\Psi,\gamma+\delta}(x)=C_{\Psi,\gamma}(x)$
this establishes the result.

\end{document}